\documentclass[11pt]{amsart} \usepackage{latexsym, amssymb, stmaryrd}

\newtheorem{thm}{Theorem}[section]
\newtheorem{lemma}[thm]{Lemma}

\newtheorem{cor}[thm]{Corollary}

\theoremstyle{definition}

\newtheorem{nrmk}[thm]{Remark}

\newtheorem{fact}[thm]{Fact}

\newcommand{\Z}{\mathbb{Z}}

\renewcommand{\o}{\Omega}

\makeatletter

\def \<{\langle}
\def \>{\rangle}

\def \z {{\mathbb Z}}
\def \*Z {{{^*}\Z}}

\def \((  {(\!(}
\def \)) {)\!)}

\def \b {\mathcal{B}}

\def \int{\operatorname{int}}
\numberwithin{equation}{section}

\def \u{\mathcal{U}}
\def \i{\iota}

\def \interior{\operatorname{interior}}

\allowdisplaybreaks[2]

\begin{document}

\title{Locally compact contractive local groups} 
\author{Lou van den Dries and Isaac Goldbring}

\address {University of Illinois at Urbana-Champaign, Department of Mathematics,
1409 W. Green St., Urbana, IL 61801}
\email{vddries@illinois.edu}

\address {University of California, Los Angeles, Department of Mathematics, 520 Portola Plaza, Box 951555, Los Angeles, CA 90095-1555, USA}
\email{isaac@math.ucla.edu}
\urladdr{www.math.ucla.edu/~isaac}

\begin{abstract}
We study locally compact contractive local groups, that is, locally compact local groups with a contractive pseudo-automorphism.  We prove that if such an object is locally connected, then it is locally isomorphic to a Lie group.  We also prove a related structure theorem for locally compact  contractive local groups
which are not necessarily locally connected.  These results are local analogues of theorems for locally compact contractive groups.
\end{abstract}
\maketitle

\begin{section}{Introduction}

\noindent
Throughout $G$ is a local group as defined in \cite{Gold}. We let $1$
be its identity, $\Lambda\subseteq G$ be the domain of its inversion map, and 
$\Omega\subseteq G\times G$ be the domain of its product map. All 
local groups in this paper are assumed to be hausdorff, and likewise,
``topological group'' means ``hausdorff topological group''.

An automorphism $\varphi$ of  
a topological group $H$ is said
to be 
\emph{contractive} if  
$$\lim_{n\to \infty} \varphi^n(x)=1\ \text{ for all }x\in H,$$ and
we call a topological group {\em contractive\/}\footnote{
\emph{Contractible} in \cite{Si}, but this term has another 
meaning in topology.}
if it has a contractive automorphism.
In \cite{Si} it is shown that locally compact connected contractive 
topological groups are (finite-dimensional, real) Lie groups.
In response to a 
question by Svetlana Selivanova we prove here a local analogue of this result. 
To formulate 
this analogue precisely, we define
 a \emph{contractive pseudo-automorphism} of $G$ to be a
morphism $\varphi\colon G\to G$ of local groups such that
for some open neighborhood $U$ of $1$ in $G$ the map $\phi|U: U \to G$ is 
injective and open, and $\lim_{n\to \infty} \varphi^n(x)=1$ for all $x\in U$. 
Call $G$ {\em contractive\/} if $G$ has a contractive
pseudo-automorphism.  

\begin{thm}\label{con} If $G$ is locally compact, 
locally connected, and contractive, then $G$ is locally isomorphic to a 
contractive Lie group.
\end{thm} 

\noindent
The recent solution \cite{Gold} of a local
version of Hilbert's 5th problem is of no help here, and we use instead
an old result due to Mal'cev~\cite{M}
to the effect that local groups satisfying a certain generalized associative law
embed into topological groups. In Section 2 we prove
Mal'cev's theorem. In Section 3 we show that if $G$ is
contractive in a strong way, then $G$ obeys 
the generalized associative law that makes Mal'cev's theorem applicable. 
In Section 4 we use this to derive 
Theorem~\ref{con} from the corresponding global result in \cite{Si}. We also
prove a related structure theorem for locally compact 
contractive local groups that are not necessarily locally connected.

\medskip\noindent
See \cite{Gold} for the definition of $G|U$ for an open neighborhood
$U$ of $1$ in $G$, and of ``morphism of local groups'' 
(also called ``local group morphism'' below). Recall also from \cite{Gold}
that two local groups are said to be locally isomorphic if they have isomorphic
restrictions to open neighborhoods of their identity. 
Here are definitions of some auxiliary
notions. Let $X\subseteq G$. We call $X$ {\em symmetric\/} if 
$X\subseteq \Lambda$ and $X^{-1}=X$; in particular, 
$G$ is symmetric iff $\Lambda=G$. The largest 
symmetric subset of $X$ is its {\em symmetrization} $X_s$: 
$$X_s:= \{x\in X\cap \Lambda:\ x^{-1}\in X\cap \Lambda\}  \qquad 
(\text{so }G_s=\Lambda\cap \Lambda^{-1}).$$
If $U$ is an open neighborhood of $1$, then so is $U_s$. If 
$\varphi: G \to G$ is a contractive pseudo-automorphism of $G$, then
$\varphi(G_s)\subseteq G_s$, and
the restriction of $\varphi$ to a
map $G_s \to G_s$ is a contractive pseudo-automorphism of $G|G_s$. 
We call $G$ {\em neat\/} if $\Lambda=G$ and $(xy, y^{-1})\in \Omega$
for all $(x,y)\in \Omega$. Note that $G|U$ is neat for any symmetric open 
neighborhood $U$ of $1$ with $U\times U\subseteq \Omega$.

\end{section}

\begin{section}{Mal'cev's theorem}

\noindent
Theorem~\ref{T:malcev} below provides a necessary and sufficient 
condition for a neat local group to admit an injective local group morphism
into a topological group.
Because some of its byproducts are useful in the next section we repeat 
Mal'cev's construction~\cite{M}, and include details omitted in Mal'cev's 
proof. Throughout we let $m,n$ range over $\mathbb{N}=\{0,1,2,\dots\}$. 

\medskip\noindent
We call $G$ \textbf{globalizable} if there is a topological group $H$ and 
an open neighborhood $U$ of the identity in $H$ such that $G=H|U$. Note that if
$G$ is globalizable and symmetric, then $G$ is neat.

\medskip\noindent
Let $a_1,\ldots,a_n,b\in G$.  We define the notion $(a_1,\ldots,a_n)\leadsto b$, by induction on $n$ as follows:
\begin{itemize}
\item If $n=0$, then $(a_1,\ldots,a_n)\leadsto b$ iff $b=1$;
\item $(a_1)\leadsto b$ iff $a_1=b$;
\item If $n>1$, then $(a_1,\ldots,a_n)\leadsto b$ iff for some $i\in \{1,\ldots,n-1\}$, there exist $b',b''\in G$ such that $(a_1,\ldots,a_i)\leadsto b'$, $(a_{i+1},\ldots,a_n)\leadsto b''$, $(b',b'')\in \o$ and $b'\cdot b''=b$.
\end{itemize}

\medskip\noindent
Informally, $(a_1,\ldots,a_n)\leadsto b$ if for some way of introducing parentheses into the sequence $(a_1,\ldots,a_n)$ all intermediate products are defined and the resulting product equals $b$.  A priori, there may be distinct $b,c\in G$ such that $(a_1,\ldots,a_n)\leadsto b$ and $(a_1,\ldots,a_n)\leadsto c$.

\medskip\noindent
We call $G$ \textbf{globally associative} if for all $a_1,\ldots,a_n,b,c\in G$ such that $(a_1,\ldots,a_n)\leadsto b$ and $(a_1,\ldots,a_n)\leadsto c$ we have $b=c$. If $G$ is globally asociative, so is its restriction $G|U$ to any open neighborhood $U$ of $1$. If there is an injective local group morphism from
$G$ into a topological group, then $G$ is globally associative. For neat $G$ 
the converse holds:

\begin{thm}\label{T:malcev} Suppose $G$ is neat and globally associative. Then
there is an injective local group morphism $\iota\colon G\to H$ into a topological group $H$ such that if $\phi\colon G\to L$ is any local group morphism into a topological group $L$, then there is a unique continuous group morphism $\tilde{\phi}\colon H\to L$ with $\tilde{\phi}\circ \iota=\phi$. 
\end{thm}

\begin{proof}
Let $G^*:=\bigcup_n G^{\times n}$ be the set of words on $G$. Consider a word 
$x=(x_1,\ldots,x_m)\in G^{\times m}$.  If $(x_i,x_{i+1})\in \o$, $1\le i<m$, 
then we call the word $$(x_1,\ldots,x_{i-1},x_ix_{i+1},x_{i+2},\ldots,x_m)\in G^{\times (m-1)}$$ a \emph{contraction of $x$ of type I}.  If also $x_{i+1}=x_i^{-1}$, then we call $$(x_1,\ldots,x_{i-1},x_{i+2},\ldots,x_m)\in G^{\times (m-2)}$$ 
a \emph{contraction of $x$ of type II}. If $(a,b)\in \o$ and $x_i=ab$, $1\le i\le m$ then  
$$(x_1,\ldots,x_{i-1},a,b,x_{i+1},\ldots,x_m)\in G^{\times (m+1)}$$
is an \emph{expansion of $x$ of type I}. Finally, for $a\in G$ and 
$0\le i \le m$ we call 
$$(x_1,\dots, x_{i},a, a^{-1}, x_{i+1},\dots, x_m)\in G^{\times(m+2)}$$   
an \emph{expansion of $x$ of type II}. Define an {\em admissible\/} 
sequence to be a finite sequence $w_1,\ldots,w_N$ of words $w_i\in G^*$ with 
$N\ge 1$ such that $w_{i+1}$ is a contraction or expansion of $w_i$, for all $i$ with $1\leq i<N$. This gives an equivalence relation 
$\sim$ on $G^*$ by: $x\sim y$ iff
there is an admissible sequence $w_1,\dots, w_N$ such that
$w_1=x$ and $w_N=y$. Let $H$ be the set of equivalence classes 
$[x]$ of elements $x=(x_1,\dots, x_m)\in G^*$. It is easy to check that we have a binary operation and a unary operation on $H$ given by 
$$[(x_1,\ldots,x_m)]\cdot [(y_1,\ldots,y_n)]:=[(x_1,\ldots,x_m,y_1,\ldots,y_n)]$$ and
$$[(x_1,\ldots,x_m)]^{-1}:=[(x_m^{-1},\ldots,x_1^{-1})].$$  Endowed with these operations, $H$ is a group with identity element $1_H=[\emptyset]$, the equivalence class of the empty sequence.  Note that also $1_H=[(1)]$.

Define $\iota\colon G\to H$ by $\iota g:=[(g)]$.  Clearly, $\iota G$ generates the group $H$.  We now show that $\iota$ is injective.  
(This is the part asserted without proof by Mal'cev \cite{M}.)  
The key to doing this is the following.

\medskip\noindent \textbf{Claim 1:}  Suppose that $x,y,z\in G^*$ and 
$x$ contracts to $y$ and $y$ expands to $z$.  Then one can also go from $x$ to 
$z$ by first expanding once or twice and then contracting once.

\medskip\noindent
There are some obvious cases where the relevant contraction and expansion 
operations ``commute'' and can just be interchanged. (This includes the case 
where $y$ is a contraction of $x$ of type II or
$z$ is an expansion of $y$ of type II.) So we can assume that
$y$ is a contraction of $x:=(x_1,\ldots,x_m)$ of type I,
$$ y=(x_1,\ldots,x_{i-1},x_ix_{i+1},x_{i+2},\ldots,x_m), \quad 1\le i<m,\ (x_i, x_{i+1})\in \Omega,$$ 
and $z$ is an expansion of $y$ of type I of the form
$$z=(x_1,\ldots,x_{i-1},a,b,x_{i+2},\ldots,x_m), \qquad (a,b)\in \o,\ ab=x_ix_{i+1}.$$ 
Now $G$ is neat, so $(ab,x_{i+1}^{-1})\in \o$ and $x_i=(ab)x_{i+1}^{-1}$.  Define 
\begin{align*} u:=&(x_1,\ldots,x_{i-1},ab,x_{i+1}^{-1},x_{i+1},\ldots,x_n), \\
v:=&(x_1,\ldots,x_{i-1},a,b,x_{i+1}^{-1},x_{i+1},\ldots,x_n).
\end{align*}  Then $u$ is an expansion of $x$ of type I, $v$ is an expansion of $u$ of type I, and $z$ is a contraction of $v$ of type II. 
This proves the claim.

\medskip\noindent
Define a {\em special\/} sequence to be an admissible sequence 
$w_1,\dots, w_N$ such that for some $M\in \{1,\dots,N\}$, $w_{i+1}$ 
is an expansion of $w_i$ for $1\le i<M$, and $w_{i+1}$ is a contraction of 
$w_i$ for $M\le i<N$.

\medskip\noindent
\textbf{Claim 2:} Let $x,y\in G^*$ and $x\sim y$. Then there is a
special sequence $w_1,\dots, w_N$ such that $w_1=x$ and $w_N=y$.

\medskip\noindent
To prove this, let  $w_1,\dots, w_n$ be any admissible sequence 
(typically, part of an admissible sequence connecting $x$ to $y$), and
suppose it is not special. Then $n\ge 3$ and we have a largest $m\in \{2,\dots, n-1\}$ such that
$w_{m-1}$ contracts to $w_m$ and $w_m$ expands to $w_{m+1}$. Apply Claim 1 to $w_{m-1}, w_m, w_{m+1}$ in the role of $x,y,z$, so $w_m$ gets replaced by 
one or two words. If the resulting admissible sequence is not yet special, 
apply the same procedure to it. We have to show that after a finite number of 
such steps we end up with a special sequence. The critical case is when
$m\in \{2,\dots,n-1\}$ is such that $w_{i}$ contracts to $w_{i+1}$ 
for $1\le i<m$, and $w_{i}$ expands to $w_{i+1}$ for $m\le i<n$.
Then the reader can easily check that after at most 
$$(n-m) + 2(n-m) + \dots + 2^{m-1}(n-m)=(2^m-1)(n-m)$$
such steps (applications of Claim 1) we obtain a special sequence.
This concludes the proof of Claim 2.

\medskip\noindent Now suppose that $a,b\in G$ and $\iota(a)=\iota(b)$, that is, $(a)\sim (b)$.  
By Claim 2, we can take $x\in G^*$ such that $x$ is obtained from $(a)$ by a finite 
succession of expansions (hence $x\leadsto a$), and $(b)$ is obtained from $x$ by a finite succession 
of contractions (hence $x\leadsto b$).  
Hence $a=b$ by global associativity, so $\iota$ is injective. Note:
$$\iota(1)=1_H,\ \i(a^{-1})=\i(a)^{-1} \text{ for }a\in G,\  \i(ab)=\i(a)\i(b) \text{ for }(a,b)\in \o.$$ 

\medskip\noindent
Let $\b$ be the set of open neighborhoods of $1$ in $G$, and $\i\b:=\{\i U \ | \ U\in \b\}$.  We verify the conditions (i)-(v) below that make $\i\b$ a neighborhood base at $1_H$ for a (necessarily unique) group topology on $H$, which by
convention includes here the requirement of being hausdorff. 

\medskip\noindent 
(i) Let $U,V\in \b$;  we need $W\in \b$ such that 
$\i W\subseteq \i U \cap \i V$.  Since $\i$ is injective, we can take $W=U\cap V$.

\medskip\noindent  
(ii) Let $U\in \b$; we need $V\in \b$ such that $\i V\cdot \i V\subseteq 
\i U$.  Choose $V\in \b$ such that $V\times V\subseteq \o$ and $V^2\subseteq U$.  Then for $g,g'\in V$, we have $$\i(g)\cdot \i(g')= \i(gg')\in \i U.$$

\medskip\noindent 
(iii) Let $U \in \b$; we need $V\in \b$ such that $(\i V)^{-1}\subseteq \i U$.  Choose $V\in \b$ such that $V^{-1}\subseteq U$, for example $V=U\cap U^{-1}$.  Then clearly $(\i V)^{-1}\subseteq \i U$.

\medskip\noindent
(iv) Let $h\in H$ and $U\in \b$;  we need $V\in \b$ such that 
$h(\i V) h^{-1}\subseteq \i U$. Since $H$ is generated by $\i G$ we can reduce to the case $h=\i g,\ g\in G$. Choose $V\in \b$ such that  $\{g\}\times V\subseteq \o$, $(gV)\times \{g^{-1}\}\subseteq \o$, and $(gV)g^{-1}\subseteq U$.

\medskip\noindent
(v) (hausdorff requirement) $\bigcap\{\i U \ | \ U\in \b\}\ =\ \{1_H\}$. This 
holds because $G$ is hausdorff.

\bigskip\noindent
With $H$ now being a topological group, $\i$ is clearly continuous at $1$.
Then the local homogeneity lemma 2.16 of \cite{Gold} yields that
$\i$ is continuous at each $a\in G$, and thus $\i$ is a local group morphism. 

\medskip\noindent
Let $L$ be any topological group and $\phi\colon G\to L$ a morphism of 
local groups. 

\medskip\noindent 
\textbf{Claim 3:} Suppose $x_1,\ldots,x_m,y_1,\ldots,y_n\in G$ and $(x_1,\ldots,x_m)\sim (y_1,\ldots,y_n)$.  Then $\phi(x_1)\cdots\phi(x_m)=\phi(y_1)\cdots\phi(y_n)$.

It is routine to verify the claim when 
$y=(y_1,\dots,y_n)$ is a contraction or expansion of $x=(x_1,\dots,x_m)$, 
and the general case then follows.

\medskip\noindent
By Claim 3 we can define a group morphism $\tilde{\phi}\colon H\to L$ by 
$$\tilde{\phi}([g_1,\ldots,g_n]):=\phi(g_1)\cdots \phi(g_n), \qquad (g_1,\dots, g_n\in G),$$ 
so $\tilde{\phi}\circ \i=\phi$. To check continuity of $\tilde{\phi}$, let 
$V$ be an open neighborhood of the identity in $L$. Then $U:=\phi^{-1}(V)$
is an open neighborhood of $1$ in $G$ and $\i U\subseteq \tilde{\phi}^{-1}(V)$, so $\tilde{\phi}^{-1}(V)$ is a neighborhood of $1_H$ in $H$. Thus $\tilde{\phi}$ is continuous. 
\end{proof}

%\noindent  In \cite{O}, a version of Mal'cev's theorem is proven, 
%again in attempt to add more detail to the proof that the map $\iota$ 
%above is injective.  However, instead of assuming that $G$ is neat, 
%it is assumed that $G$, $\Omega$, and $\Lambda$ are all connected 
%and that if $U$ is a neighborhood of $1$ in $G$, then for any $x\in G$, 
%there exists $x_1,\ldots,x_n\in U$ such that $(x_1,\ldots,x_n)\leadsto x$.  
%Incidentally, even under these assumptions, the proof in \cite{O} 
%that $\iota$ is injective is still suspect.  Also, \cite{O} further 
%assumes that $G$ has the structure of a local Lie group, allowing 
%$H$ to be endowed with the structure of a Lie group so that $\iota$ 
%provides a diffeomorphism between $G$ and $\iota G$.  
%Our construction allows for this as well.    

\medskip\noindent
Let $G$ be neat and globally associative. The universal property of 
$\iota, H$ in Theorem~\ref{T:malcev} determines $\iota, H$
up to unique isomorphism over $G$, and so, without claiming that $G$ is 
globalizable, we may call 
$H$ the \textbf{globalization} of $G$.
The construction in the
proof of the theorem and the local homogeneity lemma 2.16 of \cite{Gold}
show that $\i : G \to H$
is not just continuous but also open. 
In particular $\i G$ is open in $H$ and $\i$ is a homeomorphism onto $\i G$.
Accordingly, we identify $G$ with $\i G\subseteq H$ via $\i$.
Note that $G$ generates $H$. 
The following properties of $H$ are also evident from its construction.

\begin{lemma}

\

\begin{enumerate}
\item For any symmetric open neighborhood $U$ of $1$ in $G$ with $U\times U\subseteq \o$, we have $G|U=H|U$ $($and so $G|U$ is globalizable$)$.
\item If $G$ is connected, then $H$ is connected.
\item $G$ is locally compact if and only if $H$ is locally compact. 
\end{enumerate}
\end{lemma}

\end{section}

\noindent
{\bf Remark}. Olver~\cite{O} has another variant of Mal'cev's theorem, where 
$G$ is a local Lie group, $G$, $\Omega$, $\Lambda$ are connected,
and instead of neatness, it is assumed that for any 
$x\in G$ and neighborhood $U$ of $1$, there are 
$x_1,\ldots,x_n\in U$ such that $(x_1,\ldots,x_n)\leadsto x$.

\begin{section}{Contractive injective endomorphisms} 

\medskip\noindent
In this section $\varphi: G \to G$ is an injective morphism of local groups 
such that $\lim_{n\to \infty} \varphi^n(x)=1$ for all 
$x\in G$. (If $\varphi$ is also open, then $G$ is contractive.)

\begin{lemma}\label{L:phiparen}
Suppose that $a_1,\ldots,a_n,a\in G$ and $(a_1,\ldots,a_n)\leadsto a$.  Then  
also $(\varphi(a_1),\ldots,\varphi(a_n))\leadsto \varphi(a)$.
\end{lemma}

\begin{proof}
We proceed by induction on $n$.  The conclusion of the lemma is obvious when $n=0$ or $1$.  Suppose that $n>1$.  Choose $i\in \{1,\ldots,n-1\}$ and 
$b',b''\in G$ with $(b',b'')\in \o$ such that 
$$(a_1,\ldots,a_i)\leadsto b', \qquad (a_{i+1},\ldots,a_n)\leadsto b'', \qquad b'\cdot b''=a.$$  By the induction hypothesis we have 
$$(\varphi(a_1),\ldots,\varphi(a_i))\leadsto \varphi(b'), \qquad (\varphi(a_{i+1}),\ldots,\varphi(a_n))\leadsto \varphi(b'').$$
Also  $(\varphi(b'),\varphi(b''))\in \o$ and 
$\varphi(b')\varphi(b'')=\varphi(b'b'')=\varphi(a)$, and therefore
$(\varphi(a_1),\cdots,\varphi(a_n))\leadsto \varphi(a)$.
\end{proof}

\medskip\noindent 
The following is taken from \cite{Gold}. By recursion on $n$ we define the
relation $(a_1,\ldots,a_n)\rightarrow b$ for $a_1,\ldots,a_n,b\in G$ as follows:\begin{itemize}
\item If $n=0$, then $(a_1,\ldots,a_n)\rightarrow b$ iff $b=1$;
\item $(a_1)\rightarrow b$ iff $a_1=b$;
\item If $n>1$, then $(a_1,\ldots,a_n)\rightarrow b$ iff for all $i\in \{1,\ldots,n-1\}$ there exist $b',b''\in G$ such that $(a_1,\ldots,a_i)\rightarrow b'$, $(a_{i+1},\ldots,a_n)\rightarrow b''$, $(b',b'')\in \o$ and $b'\cdot b''=b$.
\end{itemize}

\medskip\noindent 
An easy induction on $n$ shows that for $a_1,\ldots,a_n,b,c\in G$, if 
$$(a_1,\ldots,a_n)\rightarrow b, \quad (a_1,\ldots,a_n)\leadsto c,$$ 
then $b=c$. By Lemma 2.5 of \cite{Gold} there is for each $n>0$ a
neighborhood $\u_n$ of $1$ such that for all $a_1,\ldots,a_n\in \u_n$ 
there is $b\in G$ with $(a_1,\cdots, a_n) \rightarrow b$.

\begin{cor}\label{L:globass}
$G$ is globally associative.
\end{cor}

\begin{proof}
Let $a_1,\ldots,a_n,b,c\in G$ be such that $$(a_1,\ldots,a_n)\leadsto b \text{ and }(a_1,\ldots,a_n)\leadsto c.$$  It is enough to derive
$b=c$.  By Lemma \ref{L:phiparen} we have 
$$(\varphi^m(a_1),\ldots,\varphi^m(a_n))\leadsto \varphi^m(b) \text{ and }(\varphi^m(a_1),\ldots,\varphi^m(a_n))\leadsto \varphi^m(c),$$ for all $m>0$.  
Choose $m$ so large that $\varphi^m(a_1),\ldots,\varphi^m(a_n)\in \u_n$.  It follows that $\varphi^m(b)=\varphi^m(c)$, and thus $b=c$.
\end{proof}

\noindent For the remainder of this section $L$ denotes a topological group.

\medskip\noindent
A \textbf{near-automorphism} of $L$ is an injective, continuous, open group 
morphism $L \to L$.  We call a near-automorphism $\tau\colon L\to L$ 
\textbf{contractive} if $\lim_{n\to \infty}\tau^n(x)=1$ for all $x\in L$.  

\medskip\noindent
For example,
$x\mapsto px\colon \mathbb{Z}_p\to \mathbb{Z}_p$ is a contractive near-automorphism of the compact additive group $\mathbb{Z}_p$ of $p$-adic integers, and is \emph{not} an automorphism.  Thus non-trivial compact groups may admit contractive near-automorphisms, but do not admit contractive automorphisms; 
see \cite{Si}, 1.8(b).

\begin{nrmk}\label{L:pseudolocal}
If $\tau\colon L\to L$ is a contractive near-automorphism of $L$, then $\tau$ is a contractive pseudo-automorphism of $L$ viewed as a local group.
\end{nrmk}

\begin{lemma}\label{L:connpseudo}
Suppose $\tau\colon L\to L$ is a near-automorphism.  Let $L_1$ be the connected component of $1$ in $L$.  Then $\tau(L_1)=L_1$ and so $\tau|L_1$ is an automorphism of $L_1$. If $L$ has only 
finitely many connected components, then $\tau$ is an automorphism of $L$.
\end{lemma}

\begin{proof}
Since $\tau$ is continuous and open, $\tau(L_1)$ is a connected open subgroup of $L$, and hence also closed in $L$, and thus $\tau(L_1)=L_1$. The set $L/L_1$ 
of cosets is the set of connected components of $L$. Suppose $L/L_1$ 
is finite.
Since $\tau(L_1)= L_1$, the function 
$$xL_1\mapsto \tau(x)L_1\colon  L/L_1\to L/L_1$$ is injective, hence bijective.  It follows that $\tau$ is an automorphism of $L$. 
\end{proof}

\begin{lemma}\label{L:conn}
Suppose $G$ is neat and $\varphi$ is open. Let $H$ be the globalization of $G$ and let $\tilde{\varphi}\colon H\to H$ be the unique extension of $\varphi$ to an endomorphism of $H$. Then the map $\tilde{\varphi}$ is open, and for 
$D:=\bigcup_n \ker(\tilde{\varphi}^n)$ we have:
\begin{enumerate}
\item $D$ is a discrete normal subgroup of $H$ and $\tilde{\varphi}^{-1}(D)= D$;
\item $\tilde{\varphi}$ descends to a contractive near-automorphism 
$$\varphi_D\colon H/D \to H/D, \quad \varphi_D(xD):=\tilde{\varphi}(x)D;$$ 
\item for any symmetric open neighborhood $U\subseteq G$ of $1$ with 
$U\times U\subseteq \o$, the image $\pi(U)$ of $U$ in $H/D$ is open, 
and the map
$$x\mapsto xD\ \colon\ U\to H/D$$ is an isomorphism $G|U \to
(H/D)|\pi(U)$ of local groups.
\end{enumerate}
\end{lemma}

\begin{proof} The openness of $\varphi$ gives the openness of 
$\tilde{\varphi}$.
It is easy to check that $D$ is a normal subgroup of $H$
and $\tilde{\varphi}^{-1}(D)= D$.  Each $\varphi^n$ is injective, so $D\cap G=\{1\}$, which gives (1), and so $\tilde{\varphi}$ descends to a
near-automorphism $\varphi_D$ of $H/D$. 
 To show that $\varphi_D$ is contractive, let $x\in H$ be given. Then 
$x=x_1\cdots x_m$, with $x_1,\ldots,x_m\in U$, so $\tilde{\varphi}^n(x) = \varphi^n(x_1)\cdots \varphi^n(x_m) \to 1$ as $n \to \infty$.
 Item (3) is straightforward.  
\end{proof}

\end{section}

\begin{section}{The structure of locally compact contractive local groups}

\medskip \noindent
In this section $H$ is a topological group, and $H_1$ is the connected 
component of its identity. Our aim here is to prove local analogues of the 
following two structure theorems for locally compact contractive groups.

\begin{fact}[\cite{MR}, (1.10) and \cite{Si}, Lemma 1.4] \label{L:connlie}
Each locally compact connected contractive group is a
Lie group. 
\end{fact}

\begin{fact}[\cite{Si}, Proposition 4.2]\label{L:fact}
If $H$ is locally compact and contractive, then $H$ is isomorphic as 
topological group to a product $H_1\times D$, where $D$ is a closed, totally disconnected, normal subgroup of $H$, and $H_1$ and $D$ are both contractive. 
(So $H_1$ is a Lie group.)
\end{fact}

\noindent Obviously, if a local group is locally isomorphic to a Lie group, 
then it is locally compact and locally connected. A strong converse of
this implication holds for contractive local groups:   
Theorem \ref{con} from the Introduction, which is our 
local analogue of Fact \ref{L:connlie}. To prove this converse
we need the next lemma whose proof is close to \cite{Si}, Lemma 1.4, and
whose purpose is to reduce to a situation where the results of the previous
section are applicable. 
%The purpose of the lemma is to show reduce to a situation 
%where some restric the assumption
%on $\varphi$ of the previous section is satisfied$G$ is locally compact and 
%contractive, 
%For $X\subseteq G$, we set 
%$$\varphi^{-n}(X):=\{x\in G \ | \ \varphi^n(x)\in X\},$$
%while $\varphi^n(X)$ has the usual meaning as a direct image. 
%Assuming $G$ is symmetric in the lemma is harmless by remarks at the end 
%of the introduction.
%The purpose of the lemma is to show that if $G$ is locally compact and 
%contractive, then $G$ has a neat contractive restriction.

\begin{lemma}\label{L:cpt} Suppose $G$ is locally compact, $\varphi$ is a
contractive pseudo-automorphism of $G$, and $V$ is a neighborhood of $1$ in 
$G$. Then there is an open symmetric neighborhood $U$ of $1$ in $G$ such that 
$U\subseteq V$, $U\times U\subseteq \Omega$, $\varphi(U)\subseteq U$, 
$\varphi|U: U \to G$ is open and injective, and 
$\lim_{n\to \infty} \varphi^n(x)=1$ for all $x\in U$.
\end{lemma}
\begin{proof} By restricting $G$ as indicated at the end of the Introduction we can assume that $G$ is symmetric. By shrinking $V$ we may assume in addition:
$V$ is compact, symmetric, 
$V\times V\subseteq \o$, and $V$ is
contained in an open neighborhood $W$ of $1$ in $G$ for which 
$\varphi|W: W \to G$ is open and injective, and $\lim_{n\to \infty} \varphi^n(x)=1$ for all $x\in W$. For $X\subseteq G$, we set 
$$\varphi^{-n}(X):=\{x\in G \ | \ \varphi^n(x)\in X\},$$
while $\varphi^n(X)$ has the usual meaning as a direct image; note that then 
for all $k\in \z$ we have $\varphi(\varphi^k(X))\subseteq \varphi^{k+1}(X)$. 
For $l\in \z$, set $V_l:=\bigcap_{k\le l} \varphi^k(V)$. We claim that then 
the family $(V_l)$ has the following properties:
\begin{enumerate}
\item $V_l$ is symmetric, $V_{l}\supseteq V_{l+1}$ and 
$\varphi(V_l)\subseteq V_{l+1}$;
\item $W\subseteq \bigcup_{l\in \z} V_l$, and $V_l\cap V$ has nonempty interior in $G$ for some $l$;
\item for every neighborhood $X$ of $1$ there exists $n_1\in \mathbb{N}$ such that for all $n\geq n_1$ we have
$\varphi^n(V)\subseteq X$;  
\item $(V_l \ | \ l\in \z)$ is a neighborhood base of $1$.
\end{enumerate}
 Item (1) is straightforward to check. To prove (2), let $x\in W$. We can take $n_0\in \mathbb{N}$ such that $\varphi^n(x)\in V$ for all $n\geq n_0$, so 
$x\in V_{-n_0}$. This proves $W\subseteq \bigcup_{l\in \z} V_l$. Now each $V_l$ is closed in $G$, so by Baire's theorem 
some $V_l\cap V$ has nonempty interior in $G$, which is (2).   

To prove (3), let $X$ be a neighborhood of $1$.  Take a compact symmetric 
neighborhood $A$ of $1$ with $A\subseteq V$, $VA\times A\subseteq \o$, and 
$A^2\subseteq X$. By (1) and (2), with $A$ in place of $V$, we obtain 
$n_0\in \mathbb{N}$ such that 
$$B:=\{x\in A \ | \ \varphi^n(x)\in A \text{ for all } n\ge n_0\}$$
has nonempty interior in $G$. Take $b\in \interior(B)$; so $b^{-1}\in A$.  
Let $x\in V$; then $(xb^{-1},b)\in \Omega$ and $x=(xb^{-1})b$. By the local 
homogeneity lemma 2.16 of \cite{Gold} we can take an open neighborhood 
$\u=\u_x$ of $1$ such that 
$$\{x\}\times \u, \{b\}\times \u\subseteq \o,\quad  
\{xb^{-1}\}\times b\u\subseteq \Omega,\quad  b\u\subseteq B,$$ and 
$x\u$ and $b\u$ are open neighborhoods of $x$ and $b$ respectively.  
Since $V$ is compact, we have $x_1,\ldots,x_m\in V$ such that 
$V\subseteq x_1\u_{x_1}\cup \cdots \cup x_m\u_{x_m}$. Choose $n_1\in \mathbb{N}$ such that $n_1\ge n_0$ and 
$\varphi^n(x_1b^{-1}),\ldots,\varphi^n(x_mb^{-1})\in A$ for all 
$n\geq n_1$.  Since $\varphi^n(B)\subseteq A$ for all $n\ge n_0$, and
for $i=1,\dots,m$ we have
$$  x_i\u_{x_i}=(x_ib^{-1})b\u_{x_i}, \qquad b\u_{x_i}\subseteq B,$$ 
it follows that $\varphi^n(V)\subseteq A^2\subseteq X$ for all $n\geq n_1$. 
This proves (3).

Applying (3) to $X=V$ gives $n_1\in \mathbb{N}$ such that $V\subseteq V_{-n}$
for all $n\ge n_1$, and thus $V_{-n}$ is a neighborhood of $1$ for all 
$n\ge n_1$. Since $\varphi$ is open near $1$, it follows from the last part of 
(1) that all
$V_l$ are neighborhoods of $1$. Since $V_n\subseteq \varphi^n(V)$ for all $n$,
this yields (4) as a consequence of (3).

\medskip\noindent
Thus $U:= \operatorname{interior}(V_0)$  satisfies the conclusion of
the lemma.   
\end{proof}

\medskip\noindent
\textbf{Proof of Theorem \ref{con}.} Let $G$ be locally compact and locally 
connected, and let $\varphi:G\to G$ be a contractive pseudo-automorphism. Our 
job is to show that then $G$ is locally isomorphic to a contractive Lie group. 

By Lemma~\ref{L:cpt} and a remark at the end of the Introduction we can 
reduce to the case that $G$ is neat and $\varphi$ is open and injective, with
$\lim_{n\to \infty} \varphi^n(x)=1$ for all $x\in G$.
Then by Lemma \ref{L:conn} and with $H$, $D$, $\phi_D$ as in that lemma, $G$ is locally isomorphic to the topological group $L=H/D$, which has $\phi_D$ as a
contractive near-automorphism. Since $G$ is locally connected, $G$ is then also locally isomorphic to $L_1$, the connected component of 
the identity of $L$, and $L_1$ is a contractive Lie group by Lemma \ref{L:connpseudo} and Fact \ref{L:connlie}.  \qed

\medskip\noindent To obtain the local analogue of Fact \ref{L:fact}, we need
the next lemma, which is essentially Fact \ref{L:fact} with a contractive 
near-automorphism of $H$ instead of a contractive automorphism.

\begin{lemma}\label{L:globalstructure}
Suppose $H$ is locally compact and $\tau$ is a contractive near-automorphism of $H$.  Then there exists a totally disconnected, closed, normal subgroup $P$ of $H$ such that $(h,p)\mapsto hp\colon H_1\times P \to H$ is an isomorphism of topological groups, $\tau(P)\subseteq P$, and $\tau|P$ is a contractive near-automorphism of $P$.   
\end{lemma}

\begin{proof}
By Lemma \ref{L:connpseudo}, $\tau|H_1$ is a contractive automorphism of $H_1$, so by Fact \ref{L:connlie}, $H_1$ is a Lie 
group. The remainder of the proof is just like that of Proposition 4.2 in 
\cite{Si}. 
\end{proof}

\begin{thm}\label{T:structure} Suppose $G$ is locally compact and contractive.
Then $G$ is locally isomorphic to a direct product $L\times P$, where $L$ is a 
contractive Lie group and $P$ is a totally disconnected locally 
compact group with a contractive near-automorphism.
\end{thm}

\begin{proof}
As in the proof of Theorem~\ref{con} we reduce to the case that $G$ is neat
and we have an injective open morphism $\varphi: G \to G$ of local groups
such that $\lim_{n\to \infty} \varphi^n(x)=1$ for all $x\in G$.  
Let $H$, $D$ and $\varphi_D$ be as
in Lemma~\ref{L:conn}. By that lemma, $G$ is locally isomorphic to $H/D$, and so
it remains to apply 
Lemma~\ref{L:globalstructure} to $H/D$ and
$\varphi_D$ in the role of $H$ and $\tau$.
\end{proof}

\end{section}

\bibliographystyle{plain}

\end{document}